\documentclass[10pt]{amsart}
\usepackage{amscd,amsfonts,amsthm}
\usepackage{graphicx,verbatim,mathrsfs,tikz}
\usepackage{a4wide}
\usepackage[latin1]{inputenc}
\usepackage{hyperref}
\usepackage[capitalise]{cleveref}

\newcommand{\N}{\mathbb{N}}

\newcommand{\R}{\mathbb{R}}
\newcommand{\C}{\mathbb{C}}
\newcommand{\Ac}{\mathcal{A}}

\newcommand{\Cc}{\mathcal{C}}
\newcommand{\Dc}{\mathcal{D}}

\newcommand{\Nr}{\mathscr{N}}

\newcommand{\dist}{\rho_G} 
\RequirePackage{ifthen}
\newcommand{\norm}[1]{\lVert#1\rVert}%

\newtheorem{theorem}{Theorem}
\newtheorem{proposition}[theorem]{Proposition}
\newtheorem{lemma}[theorem]{Lemma}

\newcommand{\rmi}{\mathrm{i}}
\providecommand{\union}{\cup}

\providecommand{\from}{\colon}
\providecommand{\setsize}[1]{\lvert#1\rvert}
\providecommand{\argmt}{{}\cdot{}}
\providecommand{\spr}[2]{\langle#1,#2\rangle}
\providecommand{\abs}[1]{\lvert#1\rvert}

\providecommand{\xto}{\xrightarrow}

\title[Deficiency indices for adjacency matrices]
      {Addendum: ``The problem of deficiency indices for discrete
        Schr\"odinger operators on locally finite graphs''
        [J. Math. Phys. (52), 063512 (2011)]} 

\begin{document}
\author{Sylvain Gol\'enia}
\email{sylvain.golenia@u-bordeaux1.fr}
\address{Institut de Math\'ematiques de Bordeaux, Universit\'e
Bordeaux $1$, $351,$ cours de la Lib\'eration
\\$33405$ Talence cedex, France
}
\author{Christoph Schumacher}
\email{christoph.schumacher@mathematik.tu-chemnitz.de}
\address{Fakult\"at f\"ur Mathematik der Technischen Universit\"at Chemnitz,
Reichenhainer Str.~41\\
09126 Chemnitz, Germany}
\date{\today}
\begin{abstract}
In this note we answer negatively to our conjecture concerning the
deficiency indices. More precisely, given any non-negative integer~$n$,
there is locally finite graph on which the adjency matrix has
deficiency indices $(n,n)$.
\end{abstract}

\maketitle

Given a closable and densely defined symmetric operator~$T$
acting on a complex Hilbert space, the deficiency indices of~$T$
are defined by $\eta_\pm(T):=\dim\ker(T^* \mp \rmi)\in\N\union\{+\infty\}$.
The operator~$T$ possesses a self-adjoint
extension if and only if $\eta_+(T)=\eta_-(T)$.  If this is the case,
we denote the common value by $\eta(T)$ and the self-adjoint extensions
of $T$ are parametrized by the unitary group $U(\eta(T))$, e.g.,
\cite[Section X.1]{RS}. The operator~$T$ is essentially self-adjoint,
i.e., its closure is self-adjoint, if and only if $\eta(T)=0$. In this
note we discuss the possible values of $\eta(T)$, when $T$ is the
adjacency matrix acting on a locally finite and simple graph.

We recall some standard definitions of graph theory to fix notation.
A (simple, undirected) \emph{graph} is a pair $G=(E,V)$,
where $V$ is a countable set and $E\from V\times V\to\{0,1\}$
is a symmetric function with $E(v,v)=0$ for all $v\in V$.
The elements of $V$ are called \emph{vertices}.
Two vertices $v,w\in V$ with $E(v,w)=1$ form an \emph{edge} $(v,w)$,
are \emph{neighbours}, and we write $v\sim w$.
The set of neighbours of $v\in V$ is $\Nr_G(v):=\{w\in V\mid v\sim w\}$.
The number of neighbours of~$v$ is the \emph{degree}
$d_G(v):=\setsize{\Nr_G(v)}$ of~$v$.
The graph~$G=(V,E)$ is \emph{locally finite}, if $d_G(v)<\infty$ for
all $v\in V$. 
In this note, all graphs are simple, undirected and locally finite.

A \emph{path} of length~$n\in\N$ in~$G$ is a tuple $(v_0,v_1,\dotsc,v_n)\in V^{n+1}$
such that $v_{j-1}\sim v_j$ for all $j\in\{1,\dotsc,n\}$.
Such a path \emph{connects} $v_0$~and~$v_n$ and is called
\emph{$v_0$-$v_n$-path}. 
Being connected by a path is an equivalence relation on~$V$,
and the equivalence classes are called \emph{connected components} of the graph.
A graph is \emph{connected}, if all its vertices belong to the same
connected component. 
The vertex set~$V$ of a connected graph is
equipped with the \emph{graph metric} $\dist\from V\times V\to\R$,
$\dist(v,w):=\inf\{n\in\N\mid\text{there exists a $v$-$w$-path of
  length~$n$}\}$. Note that we use the convention $0\in\N$, 
so that each vertex is connected to itself with a path of
length~$0$.

We now define trees. 
An edge $e\in V\times V$, $E(e)=1$,
in a connected graph $G=(E,V)$ is \emph{pivotal},
if the graph~$G$ with the edge~$e$ removed,
i.e.\ $(\tilde E,V)$ with $\tilde E(e)=0$
and $\tilde E(e')=E(e')$ for all $e'\in V\setminus\{e\}$, is disconnected.
A \emph{tree} is a connected graph, which has only pivotal edges.

We associate to a graph~$G$ the complex Hilbert space~$\ell^2(V)$.
We denote by $\spr\argmt\argmt$ and by $\norm{\argmt}$
the scalar product and the associated norm, respectively.
The set of complex functions with finite support in~$V$ is denoted by~$\Cc_c(G)$.
One may define different discret {operators} acting on~{$\ell^2(V)$}. 
For instance, the (physical) \emph{Laplacian} is defined by
\begin{eqnarray}\label{e:defreg}
  (\Delta_{G,\circ}f)(x):=\sum_{y\in\Nr_G(x)}\bigl(f(x)-f(y)\bigr)\text,
    \text{ with }f\in\Cc_c(G)
\end{eqnarray}
It is well known that it is symmetric and essentially self-adjoint on~$\Cc_c(G)$,
see \cite{Woj}.
\nocite{Woj2}

In this note we focus on the study of the \emph{adjacency matrix}
of~$G$, which is defined by:
\begin{eqnarray}\label{e:def}
  (\Ac_{G,\circ}f)(x):=\sum_{y\in\Nr_G(x)}f(y),\text{ with }f\in\Cc_c(G).
\end{eqnarray}
This operator is symmetric and thus closable.
We denote the closure by~$\Ac_{G}$.
We denote the domain by $\Dc(\Ac_G)$, and its adjoint by~$(\Ac_G)^*$.
Unlike the Laplacian, $\Ac$ may have several self-adjoint extensions.
We investigate its deficiency indices.  
Since the operator $\Ac_G$ commutes with complex conjugation,
its deficiency indices are equal, see \cite[Theorem X.3]{RS}. This
means that $\Ac_G$ possesses a 
self-adjoint extension. Note that $\eta(\Ac_G)=0$
if and only if~$\Ac_G$ is essentially self-adjoint on~$\Cc_c(G)$. 

In \cite{MO,Mu}, one constructs adjacency matrices
for simple trees with positive deficiency indices.
In fact, it follows from their proofs that the deficiency indices
are infinite in both references.
As a general result, a special case of
\cite[Theorem~1.1]{GS} gives that, given a locally finite simple tree~$G$,
one has the following alternative:
\begin{eqnarray}\label{e:problem}
  \eta(\Ac_G)\in\{0,+\infty\}.
\end{eqnarray}
The value of $\eta(\Ac_G)$ is discussed in \cite{GS}
and linked with the growth of the tree.

In \cite[Section~3]{MW}, one finds:
\begin{theorem}\label{t:main}
For all $n\in \N\union\{\infty\}$, there is a simple graph $G$, such that
$\eta(\Ac_G)=n$.  
\end{theorem} 
Their proof is unfortunately incomplete.
However, the statement is correct, this is aim of this note.
In \cite{MW}, they provided simple and locally finite graph $G$ such
that $\eta(\Ac_G) \geq 1$ but did not check that $\eta(\Ac_G) =1$. The
problem comes from the fact that they considered a tree. More
precisely, they refered to the works of \cite{MO,Mu}. Therefore, 
\eqref{e:problem} gives $\eta(\Ac_G)=\infty$ in their case.   
Keeping that in mind and strongly motivated by some other examples,
we had proposed a drastically different scenario
and had conjectured in \cite{GS} that that for any simple graph~$G$,
one has \eqref{e:problem}.  

We now turn to the proof of Theorem~\ref{t:main}
and therefore disprove our conjecture.
First, we show that the validity of Theorem~{t:main}
is equivalent to the existence of a simple graph~$G$ with
\begin{align}\label{e:1}
  \eta(\Ac_G)=1\text.
\end{align} 
Of course, Theorem~{t:main} in particular states the existence of~$G$.
We focus on the other implication.
We denote the positive integers with~$\N^*$.
\begin{lemma}\label{l:tensor}
Let $n\in \N^*$ and $G$ be a locally finite and connected graph. Then
there exists a locally finite and connected graph $\tilde G$ such that 
\[\eta(\Ac_{\tilde G})= n\times \eta(\Ac_G)\text.\]
\end{lemma}
\proof
Let~$\hat G:=(\hat E,\hat V)$ be the disjoint union of~$n$ copies. We have:
$\hat G:=(\hat E,\hat V)$ with $\hat V:=\{1,\dots,n\}\times V$ and 
$\hat E\bigl((i,v),(j,w)\bigr):=\delta_{i,j}E(v,w)$. Note that
$\eta(\Ac_{\hat G})= n \times \eta(\Ac_G)$ since we have a direct sum.
Take now $v_0\in
V$ and connect the copies of $G$ by adding an edge between $(i, v_0)$ and
$(i+1, v_0)$, for all $i=1, \ldots, n-1$, and  denote the resulting
graph by $\tilde G$. Note that $\Ac_{\hat G}$ is bounded perturbation
of $\Ac_{\tilde G}$. Therefore, by Proposition~{p:stab} in Appendix~\ref{s:stabind},
we have $\eta(\Ac_{\hat G})=n\times\eta(\Ac_G)$.
\qed

\begin{figure}
  \includegraphics{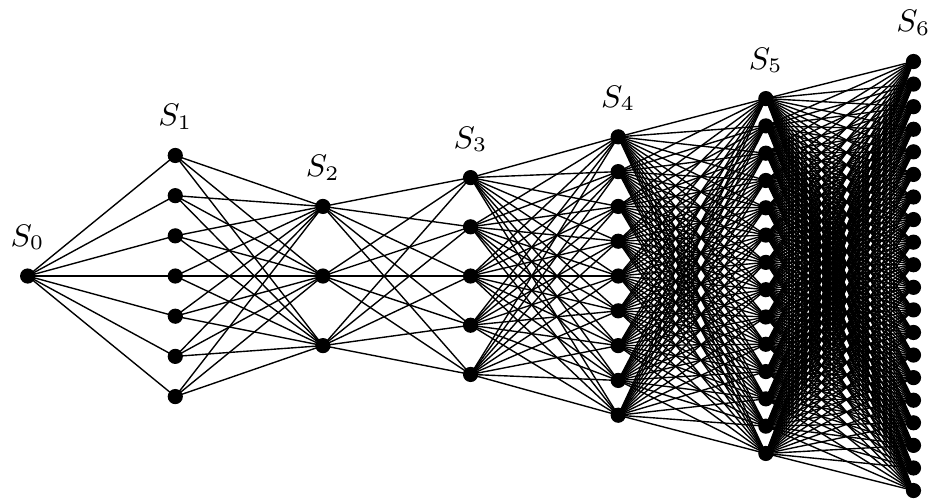}
  \caption{An antitree with spheres $S_0,\dotsc,S_6$.\label{fig_antitree}}
\end{figure}

Our example of a graph~$G$ with~\eqref{e:1}
is an antitree, a class of graphs which we define next.
See also~\cite{BK}.
The \emph{sphere} of radius~$n\in\N$ around a 
vertex~$v\in V$ is the set $S_n(v):=\{w\in V\mid d_G(v,w)=n\}$. 
A graph is an \emph{antitree}, if there exists a vertex $v\in V$
such that for all other vertices $w\in V\setminus\{v\}$
\begin{equation*}
  \Nr_G(w)=S_{n-1}(v)\union S_{n+1}(v)\text,
\end{equation*}
where $n=d_G(v,w)\ge1$.  See Figure~\ref{fig_antitree} for an example.
The distinguished vertex~$v$ is the \emph{root} of the antitree.
Antitrees are bipartite and enjoy \emph{radial symmetry},
which means that each permutation of~$V$,
which fixes the spheres around the root, induces a graph isomorphism on~$G$.

We denote the root by~$v$, the spheres by $S_n:=S_n(v)$,
{and their sizes by~$s_n:=\setsize{S_n}$.
Further, $\abs x:=\dist(v,x)$ is the distance of $x\in V$ from the root.
The operator~$P\from\ell^2(V)\to\ell^2(V)$, given by
\begin{equation*}
  Pf(x):=\frac1{s_{\abs x}}\sum_{y\in S_{\abs x}}f(y)
  \text{, for all $f\in \ell^2(V)$ and $x\in V$,}
\end{equation*}
averages a function over the spheres.
Thereby,} $P=P^2=P^*$ is the orthogonal projection
onto the space of radially symmetric functions in~$\ell^2(V)$.
A function $f\from V\to\C$ is radially symmetric, if it is constant on spheres,
i.e., for all nodes $x,y\in V$ with $\abs x=\abs y$, we have $f(x)=f(y)$.
For all radially symmetric~$f$, we define
$\tilde f\from\N\to\C$, $\tilde f(\abs x):=f(x)$, for all $x\in V$.
Note that 
\begin{equation*}
  P\ell^2(V)=\{f{\from V\to\C,\text{$f$ radially symmetric},}
    \sum_{n\in\N}s_n\abs{\tilde f(n)}^2<\infty\}
  \simeq\ell^2(\N,(s_n)_{n\in\N})\text,
\end{equation*}
where $(s_n)_{n\in\N}$ is now a sequence of weights. 
The key observation of \cite[Theorem~4.1]{BK} is that 
\begin{equation*}
  \Ac_G=P\Ac_GP
  \text{ and }
  \widetilde{\Ac_G Pf}(\abs x)
    =s_{\abs x-1}\widetilde{Pf}(\abs x-1)
    +s_{\abs x+1}\widetilde{Pf}(\abs x+1)\text,
\end{equation*}
for all $f\in \Cc_c(V)$, with the convention $s_{-1}=0$.
Using the unitary transformation
$U\from{\ell^2(\N,(s_n)_{n\in\N})}\to{\ell^2(\N)}$,
$U\tilde f(n)=\sqrt{s_n}\tilde f(n)$, we see that~$\Ac_G$
is unitarily equivalent to the direct sum of~$0$
on $(P\ell^2(V))^\perp$
and a Jacobi matrix acting on~$\ell^2(\N)$ with~$0$
on the diagonal and the sequence $(\sqrt{s_n}\sqrt{s_{n+1}})_{n\in \N}$
on the off-diagonal.

\begin{proposition}\label{p:example}
Set $\alpha>0$. Let $G$ be the  antitree with sphere sizes~$s_n$,
where $s_0:=1$, $s_n:=\lfloor{n^\alpha}\rfloor$, $n\ge1$. Then,
\begin{eqnarray*}
  \eta(\Ac_G)=
  \begin{cases}
    0,	&\text{ if }\alpha\in(0,1]\text,\\
    1,	&\text{ if }\alpha>1.
  \end{cases}
\end{eqnarray*}
\end{proposition} 
\proof Using Proposition~\ref{p:stab} from Appendix~\ref{s:stabind},
we have $\eta(\Ac_G)=\eta(J)$, where $J$ is the Jacobi matrix given by 
$a_n=\sqrt{s_ns_{n+1}}$ on the off-diagonal and $b_n=0$ on the diagonal.
Let $\tilde J$ be the Jacobi matrix given by
$\tilde a_n=\sqrt{n^\alpha(n+1)^\alpha}$ and $\tilde b_n=0$.
Now note that
\begin{align*}
  0 \le\tilde a_n-a_n&
    \le\sqrt{n^\alpha(n+1)^\alpha}-\sqrt{(n^\alpha-1)((n+1)^\alpha-1)}\\&
    =\frac{(n+1)^\alpha+n^\alpha-1}
          {\sqrt{n^\alpha(n+1)^\alpha}+\sqrt{(n^\alpha-1)((n+1)^\alpha-1)}}
    \xto{n\to\infty}1\text,
\end{align*}
therefore $\tilde a_n-a_n$ is bounded.
Hence, $\tilde J-J$ is a bounded operator, and by Proposition~\ref{p:stab},
cf.~Appendix~\ref{s:stabind}, we have $\eta(J)=\eta(\tilde J)$.

Now note $\sum_{n\in\N}\tilde a_n^{-1}=\infty$, iff $\alpha\le1$, and
\begin{align*}
  \tilde a_{n-1}\tilde a_{n+1}
    =\sqrt{(n-1)^\alpha n^\alpha}\sqrt{(n+1)^\alpha(n+2)^\alpha}&
    =\sqrt{(n^2-1)^\alpha}\sqrt{((n+1)^2-1)^\alpha}
    \le n^\alpha(n+1)^\alpha
    =\tilde a_n^2\text.
\end{align*} 
By Theorem~\ref{thm:Ber}, see Appendix~\ref{s:stabind}, applied to~$\tilde J$
we get the result.
\qed

\appendix

\section{Useful facts}
\label{s:stabind}

The theory of Jacobi matrices, as developed in \cite[Chapter VII]{Ber},
provides the following general theorem.
\begin{theorem}\label{thm:Ber}
  Let~$J$ be the Jacobi matrix with off-diagonal entries $a_n>0$
  and diagonal entries $b_n\in\R$, $n\in\N$,
  acting on $\ell^2(\N)$.
  \begin{enumerate}
    \item If $\sum_{n\in\N}a_n^{-1}=\infty$, then~$J$ is essentially
      self-adjoint on $\Cc_c(\N)$.
    \item If $\sum_{n\in\N}a_n^{-1}<\infty$, $a_{n-1}a_{n+1}\le a_n^2$
      for all $n\ge n_0$ and $\abs{b_n}\le C$
      for some constants~$n_0,C>0$,
      then~$J$ is not essentially self-adjoint on $\Cc_c(\N)$ and has
      deficiency index~$1$. 
  \end{enumerate}
\end{theorem}

We also recall that the deficiency indices are stable under the
Kato-Rellich class of perturbation and refer to
\cite[Proposition~A.1]{GS} for a proof. 

\begin{proposition}\label{p:stab}
  Given two closed and densely defined symmetric operators $S$, $T$
  acting on a complex Hilbert space and such that $\Dc(S)\subset\Dc(T)$.
  Suppose there are $a\in[0,1)$ and $b\ge0$ such that
  \begin{eqnarray}\label{e:KR}
    \norm{Tf}\le a\norm{Sf}+b\norm{f},\text{ for all }f\in\Dc(S).
  \end{eqnarray}
  Then, the closure of $(S+T)|_{\Dc(S)}$ is a symmetric operator that
  we denote by $S+T$.
  Moreover, one obtains that $\Dc(S)=\Dc(S+T)$
  and that $\eta_\pm(S)=\eta_\pm(S+T)$.
  In particular, $S+T$ is self-adjoint if and only if~$S$ is self-adjoint.
\end{proposition}

\emph{Acknowledgments:}
We would like to thank Matthias Keller for helpful discussions
and Thierry Jecko for precious remarks.

\bibliographystyle{amsalpha}
\bibliography{indices}

\providecommand{\noopsort}[1]{}\providecommand{\singleletter}[1]{#1}%
\providecommand{\bysame}{\leavevmode\hbox to3em{\hrulefill}\thinspace}
\providecommand{\MR}{\relax\ifhmode\unskip\space\fi MR }
\providecommand{\MRhref}[2]{%
  \href{http://www.ams.org/mathscinet-getitem?mr=#1}{#2}
}
\providecommand{\href}[2]{#2}
\begin{thebibliography}{MW89}

\bibitem[Ber68]{Ber}
Ju.~M. Berezanski\i, \emph{Expansions in eigenfunctions of selfadjoint
  operators}, American Mathematical Society, Providence, Rhode Island, 1968.

\bibitem[BK]{BK}
J.~Breuer and M.~Keller, \emph{Spectral analysis of certain spherically
  homogeneous graphs}, to be published.

\bibitem[GS11]{GS}
S.~Gol{\'e}nia and C.~Schumacher, \emph{The problem of deficiency indices for
  discrete schr\"odinger operators on locally finite graphs}, J. Math.\ Phys.
  \textbf{52} (2011), no.~6, 17, 063512.

\bibitem[MO85]{MO}
B.~Mohar and M.~Omladi{\v{c}}, \emph{The spectrum of infinite graphs with
  bounded vertex degrees}, Teubner, Leipzig, 1985.

\bibitem[M{\"u}l87]{Mu}
V.~M{\"u}ller, \emph{On the spectrum of an infinite graph}, Linear Algebra
  Appl. \textbf{93} (1987), 187--189.

\bibitem[MW89]{MW}
B.~Mohar and W.~Woess, \emph{A survey on spectra of infinite graphs}, J. Bull.\
  Lond.\ Math.\ Soc. \textbf{21} (1989), no.~3, 209--234.

\bibitem[RS78]{RS}
M.~Reed and B.~Simon, \emph{Methods of modern mathematical physics, tome i--iv:
  Analysis of operators}, Academic Press, 1978.

\bibitem[Woj07]{Woj}
R.~Wojciechowski, \emph{Stochastic completeness of graphs}, {Ph.D.} thesis,
  City University of New York, 2007, p.~72.

\bibitem[Woj11]{Woj2}
\bysame, \emph{Stochastically incomplete manifolds and graphs}, Progress in
  Probability \textbf{64} (2011), 163--179.

\end{thebibliography}

\end{document}